\documentclass{ws-procs9x6}

\begin{document}

\let\eps\varepsilon
\def\T{\mathbb{T}}
\def\R{\mathbb{R}}
\def\one{\mathbf{1}}

\title{Periodic homogenization with an interface}

\author{Martin Hairer and Charles Manson}

\address{Courant Institute, NYU and University of Warwick}

\begin{abstract}
We consider a diffusion process with coefficients that are periodic outside of an `interface
region' of finite thickness. The question investigated in the articles\cite{1,2} is the limiting long time / large scale
behaviour of such a process under diffusive rescaling. It is clear that outside of the interface,
the limiting process must behave like Brownian motion, with diffusion matrices given by the standard theory of homogenization.
The interesting behaviour therefore occurs on the interface.
Our main result is that the limiting process is a semimartingale whose bounded variation part is proportional to the local time spent on the interface. We also exhibit an explicit way of identifying its parameters in terms of the coefficients of the original diffusion.

Our method of proof relies on the framework provided by Freidlin and Wentzell
\cite{MR1245308} for diffusion processes on a graph in order to identify the generator of the limiting
process.
\end{abstract}

\section{Introduction}

In this note, we report on recently obtained results\cite{1,2} on the long-time large-scale behaviour of diffusions of the form
\begin{equation}\label{firstequation}
dX= b(X)\,ds + dB(s)\;,\qquad X(0) = x \in \R^d\;,
\end{equation}
where $B$ is a $d$-dimensional standard Wiener process. The drift $b$ is assumed to be smooth and such that $b(x + e_i) = b(x)$ for the unit vectors
$e_i$ with $i=2,\ldots,d$ (but not for $i = 1$). Furthermore, we assume that there exist smooth vector fields $b_\pm$ with unit period in \textit{every} direction
and $\eta > 0$ such that
\begin{equation}
b(x) = b_+(x)\;,\quad x_1 > \eta\;,\qquad b(x) = b_-(x)\;,\quad x_1 < -\eta\;.
\end{equation}
Setting $X^\eps(t) = \eps X(t/\eps^2)$, our aim is to characterise the limiting process $\bar X = \lim_{\eps \to 0} X^\eps$, if it exists.
In the sequel, we denote by
$\mathcal{L}$ the generator of $X$ and by $\mathcal{L}_\pm$ the generators of the diffusion processes $X_\pm$ given by
\eref{firstequation} with $b$ replaced by $b_\pm$. The processes $X_\pm$ will be viewed as processes on the torus $\T^{d}$,
and we denote by $\mu_{\pm}$ the corresponding invariant probability measures.
In order to obtain a diffusive behaviour for $X$ at large scales, we impose the centering condition
$\int_{\mathbb{T}^d} b_\pm(x)\,\mu_\pm(x) = 0$.

Before stating the main result, first define the various quantities involved and their relevance.

We define the `interface' of width $\eta$ by
$\mathcal{I}_{\eta} = \{x\in \mathbb{R}^d\,:\, x_{1} \in [-\eta, \eta]\}$.
In view of standard results from periodic homogenization \cite{MR2382139}, any limiting process
for $X^\varepsilon$ should behave like Brownian motion on either side of the interface $\mathcal{I}_0 = \{x_1 = 0\}$, with effective diffusion tensors given by
\begin{equation}
D^\pm_{ij} = \int_{\mathbb{T}^{d}}(\delta_{ik}+\partial_k g_i^\pm)(\delta_{kj}+\partial_k g_j^\pm)\,d\mu_{\pm}\;.
\end{equation}
(Summation of $k$ is implied.)
Here, the corrector functions $g_\pm \colon \mathbb{T}^d \to \mathbb{R}^d$ are the unique solutions to
$\mathcal{L}_\pm g_\pm = - b_\pm$, centered with respect to $\mu_\pm$. Since $b_\pm$ are centered with respect to $\mu_\pm$, such functions do indeed exist.

This justifies the introduction of a differential operator $\bar{\mathcal{L}}$ on $\mathbb{R}^d$ defined in two parts by $\bar{\mathcal{L}}_{+}$ on $I_+ = \{x_{1}>0\}$
and $\bar{\mathcal{L}}_{-}$ on $I_- = \{x_{1}<0\}$ with
\begin{equation}\label{e:defLpm}
\bar{\mathcal{L}}_{\pm} = {D_{ij}^\pm \over 2} \partial_i \partial_j\;,
\end{equation}
then one would expect any limiting process to solve a martingale problem associated to $\bar{\mathcal{L}}$. However, the above definition of $\bar{\mathcal{L}}$ is not complete, since
we did not specify any boundary condition at the interface $\mathcal{I}_0$.

In the one dimensional case \cite{1} the analysis is considerably simplified since
\begin{itemize}
\item The interface is zero dimensional in the limit and hence cannot exhibit any more complicated behavior than preferential exit behavior.
\item The non-rescaled process is time-reversible and therefore admits an invariant measure for which one has an explicit expression.
\end{itemize}
Both of these clues allow us to make a reasonable guess that in one dimension the limiting process will be some (possibly different on each side of zero) rescaling of skew Brownian motion. Since the diffusion coefficients on either side of the interface
are already determined by the theory of periodic homogenisation, the only parameter that remains to be determined is the
relative probability of excursions to either side of the interface. This can be read off the invariant measure by using the fact that the
rescaled invariant measure should converge to that of the limiting process.

One of the main ingredients in the analysis of the behavior of the limiting process at the interface is the invariant measure $\mu$ for the (original, not rescaled) process $X$. If we identify points that differ by integer multiples of $e_j$ for $j = 2,\ldots,d$, we can interpret $X$ as a process with state space $\mathbb{R} \times \mathbb{T}^{d-1}$. It then follows from the results in \cite{MR0133871} that this process admits a $\sigma$-finite invariant measure $\mu$ on $\mathbb{R} \times \mathbb{T}^{d-1}$.

Note that the invariant measure $\mu$ is \textit{not} finite and can therefore not be normalised in a canonical way. However, if we define the `unit cells' $C_j^\pm$ by
\begin{eqnarray}
C_j^+ = [j, j + 1] \times \mathbb{T}^{d-1}\\
C_j^- = [- j - 1, - j] \times \mathbb{T}^{d-1}
\end{eqnarray}
then it is possible to make sense of the quantity $q_\pm = \lim_{j \to \infty} \mu(C_j^\pm)$.

Let now $p_{\pm}$ be given by
\begin{equation}\label{e:defp}
p_{\pm}=\frac{q_\pm D^{\pm}_{11}}{q_+ D_{11}^{+}+ q_- D_{11}^{-}}\;,
\end{equation}
Unlike in the one-dimensional case, these quantities are not sufficient to characterise the limiting process since
it is possible that it picks up a non-trivial drift along the interface. It turns out that this drift can be described by drift
 coefficients $\alpha_{j}$ for $j=2,\ldots,d$ given by
\begin{equation}\label{e:valuealpha}
\alpha_j = 2 \Bigl({p_+ \over D_{11}^+} + {p_- \over D_{11}^-}\Bigr) \int_{\mathbb{R} \times \mathbb{T}^{d-1}} b_j(x)\,\mu(dx)  \;,
\end{equation}
where $\mu$ is normalised in such a way that $q_+ + q_- = 1$.

Given all of these ingredients, we construct an operator $\bar{\mathcal{L}}$ as follows. The domain
$\mathcal{D}(\bar{\mathcal{L}})$ of $\bar{\mathcal{L}}$ consists of functions $f \colon \mathbb{R}^d \to \mathbb{R}$ such that
\begin{itemize}
\item $f$ is continuous and its restrictions to $I_+$, $I_-$, and $\mathcal{I}_0$ are smooth.
\item The partial derivatives $\partial_{i} f$ are continuous for $i \ge 2$.
\item The partial derivative $\partial_{1} f(x)$ has right and left limits $\partial_{1} f|_{\pm}$ as $x \to \mathcal{I}_0$ and these limits satisfy the gluing condition
\begin{equation}
p_{+}\partial_{1} f|_{+}-p_{-}\partial_{1} f|_{-}+\sum_{j=2}^{d}\alpha_{j}\partial_{j} f=0\;. \label{derivative}
\end{equation}
\end{itemize}
For any $f \in \mathcal{D}(\bar{\mathcal{L}})$, we then set $\bar{\mathcal{L}} f(x) = \mathcal{L}_\pm f(x)$ for $x \in I_\pm$.
With these definitions at hand, we can state the main result of the article:

\begin{theorem}\label{maintheorem}
The family of processes $X^{\varepsilon}$ converges in law to the unique solution $\bar{X}$
to the martingale problem given by the operator $\bar{\mathcal{L}}$. Furthermore, there exist matrices
$M_\pm$ and a vector $K \in \mathbb{R}^d$ such that this solution solves the SDE
\begin{equation}
d\bar{X}(t) = \mathbf{1}_{\{\bar{X}_{1} \le 0\}} M_- dW(t) + \mathbf{1}_{\{\bar{X}_{1} > 0\}} M_+ dW(t) + K\,dL(t)\;.\label{e:limimt}
\end{equation}
where $L$ denotes the symmetric local time of $\bar{X}_{1}$ at the origin and $W$ is a standard $d$-dimensional
Wiener process. The matrices $M_\pm$ and the vector $K$ satisfy
\begin{equation}\label{e:parameters}
M_\pm M_\pm^T = D^\pm\;,\quad K_1 = p_+ - p_-\;,\quad K_j = \alpha_j\;,
\end{equation}
for $j =\{2,\ldots,d\}$.
\end{theorem}

\begin{figure}\label{fig:simulation}
\begin{center}
\includegraphics[height =4.7cm]{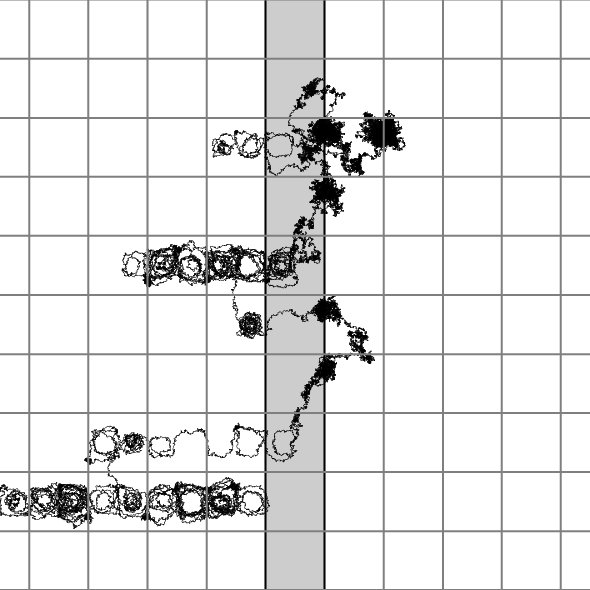}\hspace{1cm}\includegraphics[height=4.7cm]{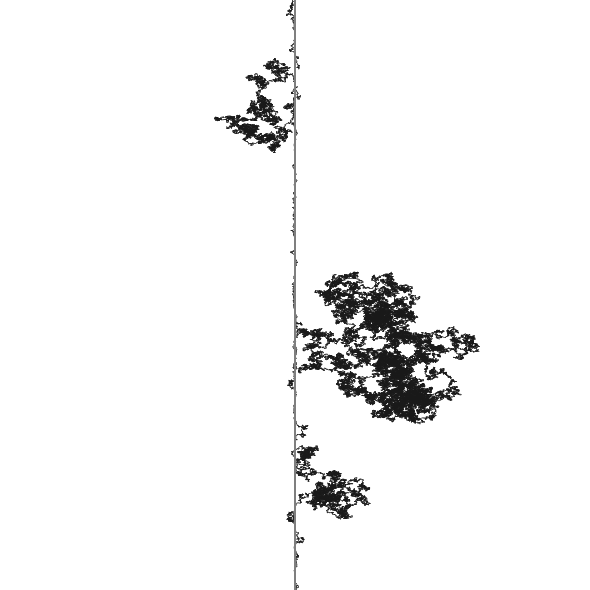}
\end{center}
\caption{Sample paths at small (left) and large (right) scales.}
\end{figure}

In Figure~\ref{fig:simulation}, we show an example of a numerical simulation of the process studied in this article.
The figure on the left shows the small-scale structure (the periodic structure of the drift is drawn as a grid). One can
clearly see the periodic structure of the sample path, especially to the left of the interface. One can also see that the effective
diffusivity is not necessarily proportional to the identity. In this case, to the left of the interface,
the process diffuses much more easily horizontally
than vertically.

The picture to the right
shows a simulation of the process at a much larger scale. We used a slightly different vector field for the drift in order
to obtain a simulation that shows clearly the strong drift experienced by the process when it hits the interface.
The remainder of this note is devoted to a short discussion of the proof of Theorem~\ref{maintheorem}.

\section{Idea of proof}
As is common in the theory of homogenization, the pattern of the proof is as follows: one first verifies tightness, then shows that any limit point satisfies the martingale problem associated to $\bar {\mathcal{L}}$, and then finally identifies solutions to this martingale problem as the unique solution to \eref{e:limimt}.

\subsection{Tightness of the rescaled processes}

We want to show that the modulus of continuity of $X^\eps$ is well-behaved uniformly in $\eps$. The only barrier to this holding can easily be shown to be the drift picked up by the process in the interface. In order to bound this, we thus need to show that the process does not spend too much time there.

We decompose the trajectory for the process $X^\varepsilon$ into excursions away from the interface,
separated by pieces of trajectory inside the interface. We first show that if the process starts inside the interface, then the expected time
spent in the interface before making a new excursion is of order $\varepsilon^2$. Then, we show that each excursion has a probability
at least $\varepsilon / \sqrt \delta$ of being of length $\delta$ or more. This shows that in the time interval $\delta$ of interest,
the process will perform at most of the order of $\sqrt \delta / \varepsilon$ excursions, so that the total time spent in the interface is indeed
of the order $\varepsilon \sqrt \delta$. Since the drift of the rescaled process is of order $1/\varepsilon$, we conclude that the
modulus of continuity is of order $\sqrt \delta$ everywhere.

\subsection{Identification of the limiting martingale problem}

In order to identify the martingale problem solved by the limiting process, it is possible to adapt
a result obtained by Freidlin and Wentzell in the context of diffusions on graphs\cite{MR1245308}.
The main ingredients are the following. For $\delta = \eps^\alpha$ with $\alpha \in ({1\over 2}, 1)$, denote by
$\tau^\delta$ the first hitting time of $\partial\mathcal{I}_\delta$ by $X^\varepsilon$. We then show that for $p_\pm$ and $\alpha_j$ as in
\eref{e:defp} and \eref{e:valuealpha}, the convergences
\begin{equation} \label{2.3}
  \mathbb{P}_{x}^{\varepsilon}[X^{\varepsilon}\bigl(\tau^{\delta}\bigr) \in I_{\pm}] \rightarrow p_{\pm}\;,\qquad
    \frac{1}{\delta}\,\mathbb{E}_{x}^{\varepsilon}\Bigl[X_{j}^{\varepsilon}\bigl(\tau^{\delta}\bigr)\Bigr]\rightarrow \alpha_{j}\;,
\end{equation}
take place  uniformly over $x\in\mathcal{I}_{\varepsilon\eta}$.

In order to show the first identity in (\ref{2.3}), let $\tau_k$ be the first hitting time of $\partial\mathcal{I}_k$ by $X$, set
$p_+^{x,k} = \mathbb{P}_x (X(\tau_k) > 0)$, and consider
\begin{equation}
\bar{p}_{+}^{k} = \sup_{x \in \mathcal{I}_\eta}p_+^{x,k}\;,\qquad\underline p_+^k = \inf_{x \in \mathcal{I}_\eta}p_+^{x,k}\;.
\end{equation}
One can then show that $\lim_{k\rightarrow\infty}|\bar{p}_{+}^{k}-\underline p_+^k|=0$ using the fact that the process returns to
any small neighborhood in $\mathcal{I}_\eta$ before $\tau_k$ with probability tending to $1$ as $k\rightarrow\infty$, allowing
the process to forget about its initial conditions through a coupling argument. The values $p_\pm$ can then be computed
in a way similar to the one-dimensional case.

The main ingredient in this calculation is the fact that the invariant measure $\mu$ for the process $X$ (which we
can view as a recurrent process on $\R \times \T^{d-1}$) gets closer and closer to multiples of $\mu_{\pm}$ away
from the interface. This can be formalised as:

\begin{proposition}\label{invariantconv}
 Let $A$ denote a bounded measurable set and denote by $\mu$ the (unique up to scaling)
 invariant $\sigma$-finite measure of the process $X$. Denote furthermore by $\mu_{\pm}$ the invariant measure of the relevant periodic process, normalised in such a way that $\mu_\pm([k,k+1]\times \mathbb{T}^{d-1}) = 1$ for
 every $k \in \mathbb{Z}$. Then there exist normalisation constants $q_\pm$ such that,
 \begin{equation}\label{e:convA}
\lim_{k \to \infty} \bigl(|\mu(A + k)-q_+\mu_{+}(A)| + |\mu(A - k)-q_-\mu_{-}(A)|\bigr) \rightarrow0\;.
 \end{equation}
(Here $k$ is an integer.) Furthermore, this convergence is exponential, and uniform over the set $A$ if we restrict its diameter.
 \end{proposition}

In order to obtain an expression for the limiting values $p_\pm$, one can now argue as follows.
Considering the first component of the limiting process, it is reasonable to expect that it converges to
a rescaling $Y$ of skew Brownian motion. This is characterised by three quantities: its diffusivity coefficients on either side of the
interface (we already know that they are given by $D_{11}^\pm$) and a parameter $p_+$ such that, setting $p_- = 1-p_+$,
\begin{equation}
\mathbb{P}_{0}^{\varepsilon}[Y\bigl(\tau^{\delta}\bigr) \in I_{\pm}] = p_{\pm}\;.
\end{equation}
The invariant measure for $Y$ is known to be proportional to Lebesgue measure on either side of the interface, with
proportionality constants $q_\pm = {p_\pm \over D_{11}^\pm}$. We can then simply solve this for $p_\pm$.

The second part of (\ref{2.3}) is shown in two steps. With $\tau_k$ as before, we have the identity
\begin{equation}\label{e:defalpha}
\alpha_j = \lim_{k \to \infty} {1\over k} \mathbb{E}_x \int_0^{\tau_k} b_j(X_s)\,ds\;,
\end{equation}
for any fixed starting point $x$ in the interface. If $k$ is large, then the process $X$ has had plenty of time to ``equilibrate'',
so that it is natural to expect that $\alpha_j$ is proportional to $\int b_j(x)\,\mu(dx)$. The only question is: what should be the correct proportionality constant?

In order to answer this question, let us assume for the sake of the argument that 
$b_j = N^{-1}\one_{[-N,N]}$ for some fixed but large value of $N$. (Note that the fact that the function $b_j$ appearing
in \eref{e:defalpha} is given by the drift of the original diffusion is irrelevant to the argument, we could ask about the value of 
this limit for any function $b$ that is localised around the interface.)
We then have
$\int b_j(x)\,\mu(dx) \approx 1$, thanks to the normalisation $q_+ + q_- = 1$. On the other hand, we know that the first component
of the rescaled process converges to skew Brownian motion described by the parameters $p_\pm$ and $D_{11}^\pm$.
Time-changing the process by a factor $D_{11}^\pm$ on either side of the origin, we can reduce ourselves to the case
of standard skew-Brownian motion with parameters $p_\pm$. Since this consists of standard Brownian motion excursions
biased to go to either side of the origin with respective probabilities $p_\pm$,
this yields in this particular example
\begin{equation}
\alpha_j = \Bigl(\frac{p_{+}}{D_{11}^{+}} + \frac{p_{-}}{D_{11}^{-}}\Bigr)\lim_{k \to \infty} {1\over k}\mathbb{E}_{0}\int_{0}^{\tau_k}\one_{[-1,1]}(B(s))\,ds\;,
\end{equation}
where $B$ is a standard Brownian motion. A simple calculation then shows that the term under the expectation
is asymptotic to $2k$, so that we do indeed recover the proportionality constant from \eref{e:valuealpha}.

\subsection{Uniqueness of the martingale problem}

Finally, in order to show uniqueness of the martingale problem, we use Theorem 4.1 from \cite{MR838085} in conjunction with the Hille-Yosida theorem to ensure that the domain of the generator to our martingale problem is large enough. It is
then possible to explicitly construct solutions to the system of SDEs given in (\ref{e:limimt}) and to show that they solve the same
martingale problem, thus concluding the proof.

\bibliographystyle{./ws-procs9x6}
\bibliography{./IMPref}
\end{document}